\documentclass{amsart}
\usepackage{amssymb,amsthm,amsmath,epsfig,latexsym,xypic}
\usepackage{calc}
\usepackage{latexsym}
\usepackage{amscd,amssymb,subfigure,hyperref}
\usepackage[arrow,matrix,graph,frame,poly,arc,tips]{xy}

\begin{document}

\newcommand{\mmbox}[1]{\mbox{${#1}$}}
\newcommand{\affine}[1]{\mmbox{{\mathbb A}^{#1}}}
\newcommand{\Ann}[1]{\mmbox{{\rm Ann}({#1})}}
\newcommand{\caps}[3]{\mmbox{{#1}_{#2} \cap \ldots \cap {#1}_{#3}}}
\newcommand{\N}{{\mathbb N}}
\newcommand{\Z}{{\mathbb Z}}
\newcommand{\R}{{\mathbb R}}
\newcommand{\Q}{{\mathbb Q}}
\newcommand{\A}{{\mathcal A}}
\newcommand{\B}{{\mathcal B}}
\newcommand{\C}{{\mathbb C}}
\newcommand{\PP}{{\mathbb P}}
\newcommand{\Fun}{\pi_1(M)}
\newcommand{\Tor}{\mathop{\rm Tor}\nolimits}
\newcommand{\ot}{\mathop{\rm OT}\nolimits}
\newcommand{\ao}{\mathop{\rm AOT}\nolimits}
\newcommand{\Ext}{\mathop{\rm Ext}\nolimits}
\newcommand{\Hom}{\mathop{\rm Hom}\nolimits}
\newcommand{\im}{\mathop{\rm Im}\nolimits}
\newcommand{\rank}{\mathop{\rm rank}\nolimits}
\newcommand{\codim}{\mathop{\rm codim}\nolimits}
\newcommand{\supp}{\mathop{\rm supp}\nolimits}
\newcommand{\CB}{Cayley-Bacharach}
\newcommand{\HF}{\mathrm{HF}}
\newcommand{\HP}{\mathrm{HP}}
\newcommand{\coker}{\mathop{\rm coker}\nolimits}
\sloppy
\newtheorem{defn0}{Definition}[section]
\newtheorem{prop0}[defn0]{Proposition}
\newtheorem{conj0}[defn0]{Conjecture}
\newtheorem{thm0}[defn0]{Theorem}
\newtheorem{lem0}[defn0]{Lemma}
\newtheorem{corollary0}[defn0]{Corollary}
\newtheorem{example0}[defn0]{Example}

\newenvironment{defn}{\begin{defn0}}{\end{defn0}}
\newenvironment{prop}{\begin{prop0}}{\end{prop0}}
\newenvironment{conj}{\begin{conj0}}{\end{conj0}}
\newenvironment{thm}{\begin{thm0}}{\end{thm0}}
\newenvironment{lem}{\begin{lem0}}{\end{lem0}}
\newenvironment{cor}{\begin{corollary0}}{\end{corollary0}}
\newenvironment{exm}{\begin{example0}\rm}{\end{example0}}

\newcommand{\msp}{\renewcommand{\arraystretch}{.5}}
\newcommand{\rsp}{\renewcommand{\arraystretch}{1}}

\newenvironment{lmatrix}{\renewcommand{\arraystretch}{.5}\small
 \begin{pmatrix}} {\end{pmatrix}\renewcommand{\arraystretch}{1}}
\newenvironment{llmatrix}{\renewcommand{\arraystretch}{.5}\scriptsize
 \begin{pmatrix}} {\end{pmatrix}\renewcommand{\arraystretch}{1}}
\newenvironment{larray}{\renewcommand{\arraystretch}{.5}\begin{array}}
 {\end{array}\renewcommand{\arraystretch}{1}}

\def \a{{\mathrel{\smash-}}{\mathrel{\mkern-8mu}}
{\mathrel{\smash-}}{\mathrel{\mkern-8mu}} {\mathrel{\smash-}}{\mathrel{\mkern-8mu}}}

\title[Resonance varieties via blowups]%
{Resonance varieties via \\blowups of $\PP^2$ and scrolls}

\author{Hal Schenck}
\thanks{Schenck supported by NSF 07--07667, NSA 904-03-1-0006}
\address{Schenck: Mathematics Department \\ University of
 Illinois \\
   Urbana \\ IL 61801\\USA}
\email{schenck@math.uiuc.edu}


\begin{abstract}
\noindent Conjectures of Suciu \cite{Su} relate the fundamental group
of an arrangement complement $M = \C^n \setminus \A$ to the 
first resonance variety of $H^*(M,\Z)$. We describe
a connection between the first resonance variety
and the Orlik-Terao algebra $C(\A)$ of the arrangement. In 
particular, we show that non-local components of $R^1(\A)$ 
give rise to determinantal syzygies of $C(\A)$. As a 
result, $Proj(C(\A))$ lies on a scroll, placing geometric
constraints on $R^1(\A)$. The key observation is that 
$C(\A)$ is the homogeneous coordinate ring associated to a 
nef but not ample divisor on the blowup of $\PP^2$
at the singular points of $\A$.
\end{abstract}
\maketitle
 
\section{Introduction}\label{sec:one}
The fundamental group of the complement $M$ of an arrangement of 
hyperplanes 
$\A = \bigcup_{i=1}^dH_i \subseteq \C^n$ 
is a much studied object. The Lefschetz-type theorem of Hamm-Le \cite{HL} 
implies that taking a generic two dimensional slice of $M$ yields an isomorphism at
the level of fundamental groups, so to study $\pi_1(M)$ we may assume
$\A \subseteq \PP^2$. Even with this simplifying assumption the situation
is nontrivial: in \cite{hirz} Hirzebruch writes ``The topology of 
the complement of an arrangement of lines in $\PP^2$ is very interesting, 
the investigation of the fundamental group very difficult''.

Presentations for $\Fun$ are given by Randell \cite{R}, 
Salvetti \cite{Sal}, Arvola \cite{Ar}, and Cohen-Suciu \cite{CS}. 
Perhaps the most compact of these is the braid monodromy 
presentation of \cite{CS}, but even this is quite complicated.
Somewhat coarser invariants of $\Fun$ are the LCS ranks and
Chen ranks. For a finitely generated group $G$, let $G=G_1$ and
define a sequence of normal subgroups inductively by $G_k = [G_{k-1},G]$.
This yields an associated Lie algebra
\[
gr(G)\otimes \Q := \bigoplus_{k=1}^{\infty} G_k/G_{k+1} \otimes \Q,
\]
with Lie bracket induced by the commutator. 
The $k$-th LCS rank $\phi_k=\phi_k(G)$ is the rank of the $k$-th quotient.
The Chen ranks of a group are the LCS ranks of the maximal metabelian
quotient $G/[[G,G],[G,G]]$. Work of Papadima and Suciu \cite{PSChen}
shows that the Chen ranks of $\Fun$ are combinatorially 
determined; but save for some special classes of arrangements, there
are no explicit formulas for either the Chen or LCS ranks. However, there are a 
beautiful pair of conjectures due to Suciu \cite{Su}, giving 
formulas for the LCS and Chen ranks in terms of the first 
resonance variety $R^1(\A)$.
The variety $R^1(\A)$ is the tangent cone at the origin to the
characteristic variety; the study of $R^1(\A)$ was pioneered by Falk
in \cite{F}. 

In the next section, we review the main subjects of investigation:
the Orlik-Solomon algebra $A=H^*(M,\Z)$, the Orlik-Terao algebra 
$C(\A)$, the first resonance variety $R^1(\A)$, and blowups 
of $\PP^2$ using certain divisors. Our main result is a description
of $C(\A)$ as the homogeneous coordinate ring of the blowup $X$ 
of $\PP^2$ at the singular points of $\A$, via a specific (nef but 
not ample) divisor $D_{\A}$. This allows us to give a geometric 
interpretation of $R^1(\A)$ in terms of certain determinantal
syzygies; we prove that if $\A$ supports a net, then $Proj(C(\A))$ lies 
on a scroll.
\section{Background}\label{sec:two}
In \cite{OS}, Orlik and Solomon gave a presentation for the cohomology ring of
the complement $M$ of a set of hyperplanes $\A \subseteq \C^n$. A 
consequence of their work is that the Betti numbers of $M$ 
are determined by the intersection lattice $L(\A)$. This lattice is 
ranked by codimension: $x \in L_i(\A)$ corresponds to a linear 
space of codimension $i$ which is an intersection of hyperplanes of $\A$.
The lattice element $\hat{0}$ corresponds to $\C^n$, and 
$y \prec x \leftrightarrow x \subsetneq y$.
We work with $\A$ central, so $\A$ defines an arrangement in both $\C^n$ and
$\PP^{n-1}$. We will depict $\A$ projectively, as below:
\begin{exm}\label{exm:braidex}
The reflecting hyperplanes of the Weyl group of $SL(4)$ are 
the six hyperplanes in $\C^4$ defined by $V(x_i-x_j)$, $1 \le i < j \le 4$. Projecting along the common 
subspace $(t,t,t,t)$ yields the {\em braid arrangement} of six planes 
containing the origin in $\C^3\!\!\!,$ or six lines in $\PP^2$:
\begin{figure}[ht]
\subfigure{%
\label{fig:braid4}%
\begin{minipage}[t]{0.35\textwidth}
\setlength{\unitlength}{16pt}
\begin{picture}(6,5)(-1,-2.5)
\put(-0.8,0){\line(1,0){5.6}}
\put(-0.8,-0.4){\line(2,1){4.7}}
\put(-0.4,-0.8){\line(1,2){2.8}}
\put(2,-0.8){\line(0,1){5.6}}
\put(4.8,-0.4){\line(-2,1){4.7}}
\put(4.4,-0.8){\line(-1,2){2.8}}
\put(-1.7,0.3){\makebox(0,0){$L_{1}$}}
\put(-1.5,-0.6){\makebox(0,0){$L_{2}$}}
\put(-0.5,-1.3){\makebox(0,0){$L_{3}$}}
\put(2,-1.3){\makebox(0,0){$L_{4}$}}
\put(4.5,-1.3){\makebox(0,0){$L_{5}$}}
\put(5.5,-0.6){\makebox(0,0){$L_{6}$}}
\end{picture}
\end{minipage}
}
\setlength{\unitlength}{0.8cm}
\subfigure{%
\label{fig:partitionlat}%
\begin{minipage}[t]{0.45\textwidth}
\begin{picture}(5,5.7)(0,-4)
\xygraph{!{0;<10mm,0mm>:<0mm,14mm>::}
[]*D(3){{\bf 0}}*-{\bullet}  
(
-@{..}[dlll]*D(3.4){123}*-{\bullet}  
(
-@{..}[d]*U(2.5){1}*-{\bullet}-@{..}[drrr]*U(2.5){\C^3}*-{\bullet} 
,-@{..}[dr]*U(2.5){2}*-{\bullet}-@{..}[drr]
,-@{..}[drr]*U(2.5){3}*-{\bullet}-@{..}[dr]
)
,-@{..}[dll]*D(2.1){25}*-{\bullet}  
(
-@{..}[d]
,-@{..}[drrrr]*U(2.5){5}*-{\bullet}-@{..}[dll]
)
,-@{..}[dl]*D(3.4){156}*-{\bullet} 
(
-@{..}[dll]
,-@{..}[drr]*U(2.5){6}*-{\bullet}-@{..}[dl]
,-@{..}[drrr]
)
,-@{..}[d]*D(2.1){36}*-{\bullet} 
(
-@{..}[dl]
,-@{..}[dr]
)
,-@{..}[dr]*D(3.4){246}*-{\bullet}  
(
-@{..}[dlll]
,-@{..}[d]
,-@{..}[drr]*U(2.5){4}*-{\bullet}-@{..}[dlll]
)
,-@{..}[drr]*D(2.1){14}*-{\bullet}  
(
-@{..}[dlllll]
,-@{..}[dr]
)
,-@{..}[drrr]*D(3.4){345}*-{\bullet}  
(
-@{..}[dllll]
,-@{..}[dl]
,-@{..}[d]
)
(
}
\end{picture}
\end{minipage}
}
\caption{\textsf{The braid 
arrangement $A_{3}$ and its intersection lattice in $\C^3$}}
\label{fig:braid}
\end{figure}
\end{exm}
\begin{defn}
The M\"{o}bius function $\mu$ : $L({\mathcal A}) \longrightarrow \Z$ is defined
 by $$\begin{array}{*{3}c}
\mu(\hat{0}) & = & 1\\
\mu(t) & = & -\!\!\sum\limits_{s \prec t}\mu(s) \mbox{, if } \hat{0} \prec t
\end{array}$$
\end{defn}
As noted, the Poincar\'e polynomial of $M$ is determined by $L(\A)$:
\[
P(M,t) = \!\!\sum\limits_{x \in L({\mathcal A})}\mu(x) \cdot (-t)^{\text{rank}(x)}.
\]
In Example~\ref{exm:braidex}, $P(M,t) = 1+6t+11t^2+6t^3$. For a central
arrangement in $\C^n$, $M \simeq \C^* \times (\PP^{n-1} \setminus \A)$, 
so by K\"unneth $P(M,t) = (1+t)P(\PP^{n-1} \setminus \A,t)$. For $n=3$, 
$b_2(M)=\sum_{p \in L_2(\A)}\mu(p)$, where $\mu(p)$ is one less than the 
number of lines through $p$.

\subsection{Orlik-Solomon algebra and $R^1(\A)$}
The Orlik and Solomon presentation for the cohomology 
ring of $M = \C^n \setminus \A$ is as follows:
\begin{defn}\label{def:OSalg}
$A=H^*(M,\Z)$ is the quotient of the exterior algebra
$E=\bigwedge (\Z^d)$ on generators $e_1, \dots , e_d$
in degree $1$ by the ideal generated by all elements of
the form $\partial e_{i_1\dots i_r}:=\sum_{q}(-1)^{q-1}e_{i_1} \cdots
\widehat{e_{i_q}}\cdots e_{i_r}$, for which
$\codim H_{i_1}\cap \cdots \cap H_{i_r} < r$.
\end{defn}
Since $A$ is a quotient of an exterior algebra, multiplication
by an element $a \in A^1$ gives a degree one differential on
$A$, yielding a cochain complex $(A,a)$:
\[
(A,a)\colon \quad
\xymatrix{
0 \ar[r] &A^0 \ar[r]^{a} & A^1
\ar[r]^{a}  & A^2 \ar[r]^{a}& \cdots \ar[r]^{a}
& A^{\ell}\ar[r] & 0}.
\]
The complex $(A,a)$ is exact as long as $\sum_{i=1}^n a_i \ne 0$;
the {\em first resonance variety} $R^1(\A)$ consists of points
$a=\sum_{i=1}^na_ie_i \leftrightarrow (a_1:\dots :a_n)$ in
$\PP(A^1) \cong \PP^{d-1}$ for which $H^1(A,a) \ne 0$.
Falk initiated the study of $R^1(A)$ in \cite{F}; among
his main innovations was the concept of a {\em neighborly
partition}: a partition $\Pi$ of $\A$ is neighborly if, for
any rank two flat $Y\in L_2(\A)$ and any block $\pi$ of $\Pi$,
\[
\mu(Y) \le |Y \cap \pi| \Longrightarrow Y\subseteq \pi, 
\]
Falk showed that all components of $R^1(A)$ arise from
such partitions, and conjectured that $R^1(\A)$ was a 
subspace arrangement. This was proved, essentially simultaneously,
by Cohen--Suciu \cite{CS} and Libgober--Yuzvinsky \cite{LY};
we will return to this in \S 4.

\subsection{The Orlik-Terao algebra}
In \cite{ot1}, Orlik and Terao introduced a commutative 
analog of the Orlik-Solomon algebra in order to answer a question of Aomoto.
\begin{defn}
Let $\A = \cup_{i=1}^d V(\alpha_i) \subseteq \PP^n$, 
and put $R=\C[y_1,\ldots,y_d]$. For each linear dependency 
$\Lambda= \sum_{j=1}^k c_{i_j}\alpha_{i_j} =0$, define 
$f_\Lambda = \sum_{j=1}^k c_{i_j} (y_{i_1}\cdots \hat
y_{i_{j}} \cdots y_{i_k})$, and let $I$ be the ideal generated 
by the $f_{\Lambda}$. The Orlik-Terao algebra $C(\A)$ is the 
quotient of $\C[y_1,\ldots,y_d]$ by $I$, and the Artinian Orlik-Terao
algebra $($the main object studied in \cite{ot1}$)$ is $C(\A)/\langle y_1^2,\ldots,y_d^2 \rangle$. 
\end{defn}
\begin{exm}\label{exm:secondex}
 Suppose $\A \subseteq {\mathbb P}^2$ is defined by 
the vanishing of $\alpha_1=x_1, \alpha_2=x_2, \alpha_3=x_3, \alpha_4=x_1+x_2+x_3$. 
The only relation is $\alpha_1+\alpha_2+\alpha_3-\alpha_4 = 0$, so
\[
C(\A) = \C[y_1,\ldots,y_4]/\langle y_2y_3y_4+y_1y_3y_4+y_1y_2y_4-y_1y_2y_3 \rangle.
\]
The homogeneous polynomial 
$y_2y_3y_4+y_1y_3y_4+y_1y_2y_4-y_1y_2y_3$ is irreducible,
hence defines a cubic surface in $\PP^3$, 
and a computation shows that the surface has four singular points. 
A classical result in algebraic geometry is that the 
linear system of four cubics through six general points in $\PP^2$ 
defines a map from the blowup of $\PP^2$ at those points to $\PP^3$ 
whose image is a smooth cubic surface. As the points move into 
special position the surface acquires 
singularities, as in this example. 
\end{exm}


In \cite{ST}, properties of the Orlik-Terao 
algebra were studied in relation to $2$--formality. An 
arrangement is $2$--formal if any dependency
among the linear forms defining the the arrangement can be obtained
as a linear combination of dependencies which involve only three
of the forms. Among the classes of $2$--formal arrangements are 
$K(\pi,1)$ arrangements and free arrangements.
However, an example of Yuzvinsky \cite{yu1} 
shows that $2$--formality is not determined by the intersection 
lattice $L(\A)$. The main result of \cite{ST} is that 
$2$--formality is determined by the quadratic component of the
Orlik-Terao ideal; the key is a computation on the tangent
space of $V(I_2) \cap (\C^*)^{d-1}.$

\subsection{Blowups of $\PP^2$}
Fix points  $p_1, \ldots p_n \in \PP^2$, and let 
\begin{equation}\label{BlowP2}
X \stackrel{\pi}{\longrightarrow}\PP^2
\end{equation} 
be the blow up of ${\PP}^2$ at these points. Then $Pic(X)$ is
generated by the exceptional curves $E_i$ over the points $p_i$, 
and the proper transform $E_0$ of a line in ${\PP}^2$. 
A classical geometric problem asks for a relationship between 
numerical properties of a divisor $D_m = mE_0-\sum a_iE_i$ on $X$, 
and the geometry of $$X \stackrel{\phi}{\longrightarrow} \PP(H^0(D_m)^\vee).$$
First, some basics. Let $m$ and $a_i$ be non-negative, 
let $I_{p_i}$ denote the ideal of a point $p_i$, and define 
\begin{equation}\label{idealJ}
J = \bigcap\limits_{i=1}^n I_{p_i}^{a_i} \subseteq \C[x,y,z]=S.
\end{equation}
Then $H^0(D_m)$ is isomorphic to the $m^{th}$ graded piece $J_m$ 
of $J$ (see \cite{h}). In \cite{DG}, Davis and 
Geramita show that if 
$\gamma(J)$ denotes the smallest degree $t$ such that $J_t$
defines $J$ scheme theoretically, then $D_m$ is very ample 
if $m > \gamma(J)$, and if $m=\gamma(J)$, then $D_m$ is very ample 
iff $J$ does not contain $m$ collinear points, counted with multiplicity. 
Note that $\gamma(J) \le reg(J)$. Now suppose that 
$\A=\cup_{i=1}^dL_i \subseteq \PP^2$, and fix 
defining linear forms $\alpha_i$ so that $L_i=V(\alpha_i)$.
Let $X$ denote the blowup of $\PP^2$ at $Sing(\A) = L_2(\A)$. The central object of our investigations is the divisor 
\begin{equation}\label{DA}
D_{\A} = (d-1)E_0 - \!\!\sum\limits_{p_i \in L_2(\A)} \mu(p_i)E_i.
\end{equation}
\subsection{Main results}
For an arrangement $\A \subseteq \PP^2$, let 
\begin{equation}\label{ProjD}
X \stackrel{\phi_{\A}}{\longrightarrow} \PP(H^0(D_{\A})^\vee).
\end{equation} 
We show that $C(\A)$ is the homogeneous coordinate ring of $\phi_{\A}(X)$, and that $\phi_{\A}$
is an isomorphism on $\pi^*(\PP^2 \setminus \A$), contracts the lines of $\A$ to
points, and blows up the singularities of $\A$. Combining results
of Proudfoot-Speyer \cite{ps} and Terao \cite{ter}, we bound the
Castelnuovo-Mumford regularity of $C(\A)$. Finally, we interpret the 
resonance varieties studied in \cite{CS}, \cite{F}, \cite{FY}, \cite{LY}, \cite{SS}, \cite{yu2} in terms of linear subsystems of $D_{\A}$, and 
connect these jump loci to linear syzygies on $C(\A)$.

\section{Connecting $H^0(D_{\A})$ to the Orlik-Terao algebra}\label{sec:three}
Let $\alpha = \prod_{i=1}^d \alpha_i$ and define a map
$R=\C[y_1,\ldots,y_d] \longrightarrow \C[ 1/ \alpha_1,\ldots,  1/ \alpha_d] = T$.
The kernel of this map is the OT ideal (see \cite{ST}), so 
$C(\A) \simeq T$. In \cite{ter}, Terao proved that the 
Hilbert series for $T$ is
given by 
\begin{equation}\label{TeraoE}
HS(T,t) = P\Big(\A, \frac{t}{1-t}\Big).
\end{equation}
In this section, we show that for $n=2$, $C(\A)$ is the 
homogeneous coordinate ring of the image of $X \stackrel{\phi_{\A}}{\longrightarrow} \PP(H^0(D_{\A})^\vee)$, with $X$ as in Equation~\ref{BlowP2}. 
For brevity, let $l_i = \alpha / \alpha_i$.
\begin{lem}\label{schemeT}
The ideal $L = \langle l_1,\ldots,l_d \rangle$ defines 
\[
\bigcap\limits_{p_i \in L_2(\A)} I_{p_i}^{\mu(p_i)} \mbox{ scheme-theoretically}.
\] 
\end{lem}
\begin{proof}
Localize at $I_p$, where $p \in L_2(\A)$. Then in $S_{I_p}$, $\alpha_i$ is a unit
if $p \not \in V(\alpha_i)$. Without loss of generality, suppose forms $\alpha_1, \ldots 
\alpha_m$ vanish on $p$, and the remaining forms do not. Thus,
\[
L_{I_p} = \langle \alpha_2\cdots\alpha_m, \alpha_1\cdot \alpha_3\cdots \alpha_m, \ldots,  \alpha_1\cdots \alpha_{m-1} \rangle.
\]
Now note that $I_p^{\mu(p)}$ has $\mu(p)+1$ generators of degree $\mu(p)$. 
Since $\mu(p) = m-1$ and the forms in $L_{I_p}$ are linearly independent, 
equality follows.
\end{proof}
\begin{lem}\label{idealT}
The minimal free resolution of $S/L$ is
\begin{small}
\[
0 \longrightarrow S(-d)^{d-1} \xrightarrow{\psi}
S(-d+1)^d
\xrightarrow{\left[ \!\begin{array}{ccc}
l_1,& \cdots& ,l_d
\end{array}\! \right]}
 S \longrightarrow S/L  \longrightarrow 0, \mbox{ where }
\]
\end{small}
\begin{tiny}
\[
\psi = {\left[ \!
\begin{array}{ccccc}
\alpha_1  & 0          & \cdots & \cdots & 0        \\
-\alpha_2 & \alpha_2   & 0      & \cdots  & \vdots   \\
0         &  -\alpha_3 & \ddots & \ddots  & \vdots         \\ 
\vdots    & 0          &        & \ddots        &  0       \\ 
\vdots    & \vdots     &  \ddots      &         & \alpha_{d-1}  \\ 
0         & \cdots     & \cdots &  0      & -\alpha_d
\end{array}\! \right]}
\]
\end{tiny}
\end{lem}
\begin{proof}
The columns of $\psi$ are syzygies on $L$. Since 
the maximal minors of $\psi$ generate $L$, the 
result follows from the Hilbert-Burch theorem and Lemma~\ref{schemeT}.
\end{proof}
\vskip .1in
\begin{thm}~\label{sectionsD} 
$H^0(D_{\A}) \simeq Span_{\C}\{l_1,\ldots, l_{d} \}$ and 
$H^1(D_{\A})=0=H^2(D_{\A}).$
\end{thm}
\begin{proof}
The remark following Equation~\ref{idealJ} shows 
that $H^0(D_{\A}) \simeq J_{d-1}$. Since $K = -3E_0+\sum E_i$,
by Serre duality 
\[
H^2(D_{\A})\simeq H^0((-d-2)E_0+\!\!\!\!\sum\limits_{p_i \in L_2(\A)} (\mu(p_i)+1)E_i),
\]
which is clearly zero. Using that $X$ is rational, it follows 
from Riemann-Roch that 
$$h^0(D_{\A})-h^1(D_{\A})=\frac{D_{\A}^2-D_{\A} \cdot K}{2}+ 1.$$
The intersection pairing on $X$ is given by 
$E_i^2  =  1$ if $i=0$, and $-1$ if $i \ne 0$, and 
\[
E_i \cdot E_j  =  0 \mbox{ if } i \ne j.
\]
Thus,
\begin{equation}
\begin{aligned}
D_{\A}^2 &= (d-1)^2 - \!\!\sum\limits_{p \in L_2(\A)} \!\!\!\mu(p)^2\\
-D_{\A}K    & = 3(d-1) - \!\!\sum\limits_{p \in L_2(\A)} \!\!\!\mu(p),
\end{aligned}
\end{equation}
yielding
\begin{equation}
\begin{aligned}
h^0(D_{\A})-h^1(D_{\A}) &= \frac{\Big(d-1\Big)^2-\sum \mu\Big(p\Big)^2 + 3\Big(d-1\Big) - \sum \mu\Big(p\Big)}{2}+ 1 \\
& =  {d+1 \choose 2} -\!\!\sum\limits_{p \in L_2(\A)} \!\!\!{\mu(p)+1 \choose 2}.
\end{aligned}
\end{equation}
Double counting the edges between $L_1(\A)$ and $L_2(\A)$ yields
$$
{d \choose 2} = \sum_{p \in L_2(\A)} \!\!\! {\mu(p)+1 \choose 2},
$$
hence $h^0(D_{\A})-h^1(D_{\A}) = d$. From 
Lemmas~\ref{schemeT} and~\ref{idealT} the Hilbert function satisfies
\[
d = HF(\langle l_1,\ldots,l_d \rangle, d-1) = HF(\bigcap\limits_{p_i \in L_2(\A)} I_{p_i}^{\mu(p_i)}, d-1).
\] 
The observation after Equation~\ref{idealJ} now implies that $h^0(D_{\A})=d$.
\end{proof}
\noindent It follows from Theorem~\ref{sectionsD} that $C(\A)$ is 
the coordinate ring of $\phi_{\A}(X)$. Note also that by 
Lemmas~\ref{schemeT} and \ref{idealT}, the constant $\gamma(L) = d-1$, so 
$$dE_0 - \!\!\sum\limits_{p_i \in L_2(\A)} \mu(p_i)E_i$$
is very ample, and gives a De Concini-Procesi wonderful model \cite{DP}:
a compactification $\overline{M}$ of $M$ such that $\overline{M} \setminus M$ 
is a normal crossing divisor. 
However, since every line of $\A$ contains exactly $d-1$ 
points counted with multiplicity, the divisor 
$D_{\A}$ is not very ample. The description of $C(\A)$ makes it 
obvious that $V(I)$ is an irreducible, 
nondegenerate rational variety, and by \cite{ST} 
$V(I) \setminus V(y_1\cdots y_d)$ is smooth. 
Here is a more explicit description of the map:
\begin{thm}\label{geometrySurf}
The map $\phi_{\A}$ 
\begin{enumerate}
\item is an isomorphism on $\pi^*(\PP^2 \setminus \A)$.
\item contracts the lines of $\A$ to points on $X$.
\item takes $E_p$ to a rational normal curve of degree $\mu(p)$.
\end{enumerate} 
\end{thm}
\begin{proof}
For the first part, without loss of generality suppose that $\alpha_1 \cdot \alpha_2 \cdot \alpha_3 = xyz$
and write $L = \prod_{i=4}^d \alpha_i$. Then 
\[
\phi_{\A}  = [yzL, xzL, xyL, \ldots ].
\]
Thus, the first three entries of $\phi_{\A}$ define the Cremona transformation, which is an isomorphism from
$\PP^2 \setminus V(xyz)$ to itself. Since $\PP^2 \setminus \A $ is contained in $\PP^2 \setminus V(xyz)$, $(1)$ follows. For $(2)$, 
suppose $p$ is a point of $V(\alpha_i)$. Since $\alpha_i$ divides 
$l_j$ for all $j \ne i$, this means $l_j(p) = 0$ if $j \ne i$. 
Hence $\phi_{\A}(V(\alpha_i))$ is the $i^{th}$ coordinate point of $\PP^{d-1}$. The final part
follows from the fact that $D_{\A}|_{E_p}$ is a divisor on $E_p$ of degree $D_{\A} \cdot E_p = \mu(p)$, and
$E_p \simeq \PP^1$.
\end{proof}
\subsection{Castelnuovo--Mumford regularity and graded betti numbers}
The Castelnuovo--Mumford regularity of
a coherent sheaf $\mathcal{F}$ on $\PP^n$ is usually phrased in terms
of vanishing of certain cohomology modules. Letting 
$N = \oplus_nH^0(\mathcal{F}(n))$, we may \cite{eis} rephrase
the condition as
\begin{defn}
For a polynomial ring $R$, a finitely generated, graded $R$--module $N$ has
Castelnuovo-Mumford regularity $j$ if $j$ is the
smallest number such that $Tor_i^R(N,\C)_{i+j+1} = 0$ for all $i$. 
The graded betti numbers of a graded $R$--module $N$ are indexed by 
\[
b_{ij} = \dim_{\mathbb{C}}\Tor_i^R(N,\C)_{j}.
\]
\end{defn}
\begin{exm}\label{braidex1}
We revisit Example~\ref{exm:braidex}. The four triple points yield four 
quadratic generators for the Orlik-Terao ideal $I$. These four
quadrics generate $I$ (see \cite{ST}), and a 
computation in {\tt Macaulay2} yields the  
graded betti numbers of $C(\A)$:
\begin{small}
$$
\vbox{\offinterlineskip 
\halign{\strut\hfil# \ \vrule\quad&# \ &# \ &# \ &# \ &# \ &# \ &# \ &# \ &# \ &# \ &# \ &#  \cr
total&1&4&5&2\cr \noalign {\hrule} 0&1 &--&--& --&-- \cr 1&--&4 &2 & --&--
\cr 2&--&--&3 &2 \cr \noalign{\smallskip} }}
$$
\end{small}
\noindent This diagram is read as follows: the entry in position $(i,j)$ is simply $b_{i,i+j}$, e.g. 
\[
\dim_{\mathbb{K}}\Tor_2^R(C(\A),\C)_4 = 3.
\]
The betti table has a very nice interpretation in terms of 
Castelnuovo-Mumford regularity: the regularity is the index of
the last nonzero row. 
\end{exm}
\begin{thm}\label{reg}
For $\A \subseteq \PP^n$, $C(\A)$ is $n$--regular.
\end{thm}
\begin{proof}
In \cite{ps}, Proudfoot and Speyer show that the
Orlik-Terao algebra is Cohen-Macaulay (for $n=2$ this also follows from 
Theorem~\ref{sectionsD}). Thus, there exists a regular sequence on $C(\A)$ of 
$\dim(V(I))+1 = n+1$ linear forms; quotienting by this
sequence yields an Artinian ring whose Hilbert series is the
numerator of the Hilbert series of $C(\A)$. The 
regularity of an Artinian module is equal to the length
of the module, so the result follows from Equation~\ref{TeraoE}.
\end{proof}
It follows easily from Theorem~\ref{reg} and Terao's work in \cite{ter} that
\begin{prop}\label{b23}
For $\A \subseteq \PP^n$ with $|\A| = d$, if $I=I_2$, then 
 \[
\dim_{\C}Tor_2^R(C(\A),\C)_{3} = 2 \Big( {d \choose 3}-1 \Big) - \Big( d-3 \Big)\Big(\!\!\sum\limits
_{p_i \in L_2(\A)}\!\! \mu(p_i)+1 \Big).
\]
\end{prop}
\begin{exm}\label{braidex2}
The $A_3$ arrangement is supersolvable, so by \cite{ST} 
$I=I_2$, and Proposition~\ref{b23} shows there are two
linear first syzygies on $I$. This explains the top row of the betti
table in Example~\ref{braidex1}.
\end{exm}
\section{Nets, syzygies, and scrolls}
In \cite{LY}, Libgober-Yuzvinsky found a surprising connection
between nets and the first resonance variety. The approach was 
further developed by Yuzvinsky in \cite{yu2}, with a beautiful 
complete picture emerging in Falk and Yuzvinsky's paper on 
multinets \cite{FY}. In this section, we connect nets to 
the linear syzygies of $C(\A)$, and hence to $R^1(\A)$. This
allows us to give an interpretation of the first resonance 
variety in terms of the geometry of $X_A$.

Suppose $Z$ is a subset of the intersection points of $\A$, and 
let $J$ denote the 
$|Z| \times d$ incidence matrix of points and lines and 
$E$ denote a $d \times d$ matrix with every entry one. 
If $\widehat{Z}$ is the blowup of $\PP^2$ at the points of $Z$,
then \cite{LY} shows that 
\[
J^tJ - E = Q(\widehat{Z})
\] 
is the intersection form on $\widehat{Z}$, and is a generalized Cartan matrix.
Using the Vinberg classification of such matrices \cite{kac}, they show that any
component of $R^1(A)$ corresponds to a choice of points $Z$ such that
$Q(\widehat{Z})$ consists of at least three affine blocks, with no finite or
indefinite blocks, and the block sum decomposition of $Q(\widehat{Z})$ yields
a neighborly partition. 
Before going into the details of the connection between multinets,
divisors and syzygies, we give a pair of motivating examples.
\begin{exm}\label{exm:p2}
The matroid $(9_3)_2$ of Hilbert and Cohn-Vossen is realized
below by $\A = V(xyz(x+y)(y+z)(x+3z)(x+2y+z)(x+2y+3z)(2x+3y+3z))$.
It has nine triple points and nine double points, thus 
$P(\A,t)=(1+t)(1+8t+19t^2)$.
\begin{figure}[ht]
\subfigure{%
\label{fig:dnonpap}%
\begin{minipage}[t]{0.3\textwidth}
\setlength{\unitlength}{9pt}
\begin{picture}(5,7.7)(-2,-1.35)
\multiput(-1,1.5)(0,1){2}{\line(1,0){9.2}}
\multiput(1,-1.5)(3,0){2}{\line(0,1){8}}
\put(-1,3.5){\line(2,-1){9.2}}
\put(-1,4.5){\line(2,-1){9.2}}
\put(-0.95,4.8){\line(3,-2){9.2}}
\put(0,6.5){\line(1,-1){8}}
\put(1,-2.2){\makebox(0,0){$1$}}
\put(4,-2.2){\makebox(0,0){$2$}}
\put(8,-2.2){\makebox(0,0){$3$}}
\put(8.85,-2.1){\makebox(0,0){$4$}}
\put(9,-1.2){\makebox(0,0){$5$}}
\put(9,-0.3){\makebox(0,0){$6$}}
\put(9,1.5){\makebox(0,0){$7$}}
\put(9,2.5){\makebox(0,0){$8$}}
\end{picture}
\end{minipage}
}
\caption{\textsf{An arrangement realizing $(9_3)_2$}}
\label{fig:nonpap}
\end{figure}

The graded betti numbers for $C((9_3)_2)$ are:

\begin{small}
$$
\vbox{\offinterlineskip 
\halign{\strut\hfil# \ \vrule\quad&# \ &# \ &# \ &# \ &# \ &# \ &# \ &# \ &# \ &# \ &# \ &# \ &# \ &#\ \cr
total&1&11&75&156&145&66&12 \cr \noalign {\hrule} 0&1 &--&--& --&--&-- &-- \cr 1&--&9 &-- & --&--&--&--   \cr 2&--&2 &75 &156 &145&66&12  \cr \noalign{\bigskip} \noalign{\smallskip} }}
$$
\end{small}
\end{exm}
\begin{exm}\label{exm:p1}
The $(9_3)_1$ matroid of Hilbert and Cohn-Vossen is realized
below by $\A = V(xyz(x-y)(y-z)(x-y-z)(2x+y+z)(2x+y-z)(2x-5y+z))$. 
It has nine triple points and nine double points, so $P((9_3)_1,t) = P((9_3)_2,t)$.
\begin{figure}[ht]
\subfigure{%
\label{fig:dpap}%
\begin{minipage}[t]{0.3\textwidth}
\setlength{\unitlength}{8.5pt}
\begin{picture}(5,9)(-2,-1)
\multiput(-1.5,1)(0,2){2}{\line(1,0){9}}
\multiput(1,-1)(2,0){2}{\line(0,1){10}}
\put(-1,7){\line(3,-4){6}}
\put(0,9){\line(3,-4){7}}
\put(0,-1){\line(1,2){4.5}}
\put(-1.5,6){\line(3,-2){9}}
\put(-0.35,-1.8){\makebox(0,0){$1$}}
\put(1,-1.8){\makebox(0,0){$2$}}
\put(3,-1.8){\makebox(0,0){$3$}}
\put(5.15,-1.8){\makebox(0,0){$4$}}
\put(7.2,-1.1){\makebox(0,0){$5$}}
\put(8.05,-0.4){\makebox(0,0){$6$}}
\put(8.15,1){\makebox(0,0){$7$}}
\put(8.15,3){\makebox(0,0){$8$}}
\end{picture}
\end{minipage}
}
\caption{\textsf{An arrangement realizing $(9_3)_1$}}
\label{fig:pap}
\end{figure}

However, the graded betti numbers for $C((9_3)_1)$ are:

\begin{small}
$$
\vbox{\offinterlineskip 
\halign{\strut\hfil# \ \vrule\quad&# \ &# \ &# \ &# \ &# \ &# \ &# \ &# \ &# \ &# \ &# \ &# \ &# \ &#\ \cr
total&1&13&77&156&145&66&12 \cr \noalign {\hrule} 0&1 &--&--& --&--&-- &-- \cr 1&--&9 &2 & --&--&--&--   \cr 2&--&4 &75 &156 &145&66&12  \cr \noalign{\bigskip} \noalign{\smallskip} }}
$$
\end{small}
The arrangement $(9_3)_1$ possesses a pair of linear first syzygies, while
$(9_3)_2$ has no linear first syzygies. An easy check shows that  $(9_3)_1$ 
admits a neighborly partition $|169|258|347|$ and has a corresponding
non-local component (see below) in $R^1(A)$, whereas $(9_3)_2$ does not. 
To better understand the connection between $R^1(A)$ and syzygies, 
we now review two constructions.
\end{exm}
\subsection{Nets and multinets}
It is easy to see that any $p \in L_2(\A)$ with $\mu(p) \ge 2$ yields a 
component $\PP^{\mu(p)-1} \subseteq R^1(A)$. Such components are 
called {\em local components}. Components which are not of this 
type are called {\em essential}. In \cite{yu2}, 
Yuzvinsky used {\em nets} to analyze the essential components of $R^1(A)$.
\begin{defn}~\label{Net}
Let $3 \le k \in \Z$. A {\it $k$-net} in $\PP^2$ is a pair $(\A,Z)$
where $\A$ is a finite set of distinct lines partitioned into $k$ subsets 
$\A=\bigcup_{i=1}^k\A_i$ and $Z$ is a finite set of points, such that:
\begin{enumerate}
\item{for every $i \ne j$ and every $L \in\A_i, \ L'\in\A_j$, $L \cap
L'\in Z$.}
\item{for every $p \in Z$ and every $i \in \{1,\ldots,k\}$, $\exists$ a unique
$L\in A_i$ containing $Z$}.
\end{enumerate}
\end{defn}
Thus, for a $k$-net, $|A_i|=|L \cap Z|$ for any block $A_i$ and line $L \in \A$;
denote this number by $m$. Following Yuzvinsky, we call $m$ the order of the net,
and refer to a $k$-net of order $m$ as a $(k,m)$-net; note that $|Z|=m^2$. 
Yuzvinsky shows in \cite{yu2} that a net must have $k \in \{3,4,5\}$, 
and improves this in \cite{yu3} to $k \in \{3,4\}$. 

In \cite{FY}, Falk and
Yuzvinsky extend the notion of a net to a multinet; in a multinet 
lines may occur with multiplicity. Write $\A_w$ for a multiarrangement, where 
$w \in \N^{d}$, and $w(L)$ denotes the multiplicity of a line. 
\begin{defn}~\label{MNet}
A weak $(k,m)$-multinet on a multi-arrangement $\A_w$ is a pair $(\Pi,Z)$ where $\Pi$ 
is a partition of $\A_w$  into $k\geq 3$ classes $A_1, \ldots ,A_k,$ and $Z$ is a set of 
multiple points, such that
\begin{enumerate}
\item $\sum_{L \in A_i} w(L)=m,$ independent of $i.$
\item For every $L \in A_i$ and $L'\in A_j,$ with $i\ne j$, $L \cap L' \in Z.$ 
\item For each $p\in Z$, $\sum_{L \in A_i, p \in L} w(L)$ is a constant $n_p$, independent of $i$.
\end{enumerate}
A multinet is a weak multinet satisfying the additional property
\begin{enumerate}
\setcounter{enumi}{3}
\item For $i \in \{1,\ldots, k \}$ and $L,L'\in A_i,\mbox{ }\exists$ a sequence $L=L_0,L_1,\ldots,L_r=L'$ such that $L_{j-1}\cap L_j\not \in Z$ for $1 \leq j \leq r.$
\end{enumerate}
\end{defn}
\begin{exm}\label{exm:b3}
The reflection arrangement of type $B_3$ is depicted 
below (there is also a line at infinity). Falk and Yuzvinsky show 
that this arrangement supports a multinet which is not a net:
assign weight two to lines $(3,6,8)$ and weight one
to the remaining lines.
\begin{figure}[ht]
\subfigure{%
\label{fig:B3-a}%
\begin{minipage}[t]{0.3\textwidth}
\setlength{\unitlength}{0.58cm}
\begin{picture}(5,4.8)(-0.2,-1)
\multiput(1,0)(1,0){2}{\line(1,1){3}}
\multiput(4,0)(1,0){2}{\line(-1,1){3}}
\multiput(2.5,0)(0.5,0){3}{\line(0,1){3}}
\put(1,1.5){\line(1,0){4}}
\put(4.2,-0.5){\makebox(0,0){$1$}}
\put(5.2,-0.5){\makebox(0,0){$2$}}
\put(5.5,1.5){\makebox(0,0){$3$}}
\put(5.2,3.5){\makebox(0,0){$4$}}
\put(4.2,3.5){\makebox(0,0){$5$}}
\put(3.5,3.5){\makebox(0,0){$6$}}
\put(3,3.5){\makebox(0,0){$7$}}
\put(2.5,3.5){\makebox(0,0){$8$}}
\end{picture}
\end{minipage}
}
\caption{\textsf{The $\operatorname{B}_3$-arrangement}}
\label{fig:B3arr}
\end{figure}
\end{exm}
The following lemma of \cite{FY} will be useful:
\begin{lem} 
\label{numerology}
Suppose $(\A_w,Z)$ is a weak $(k,m)$-multinet. Then
\begin{enumerate} 
\item $\sum_{L\in \A_w} w(L)=km.$
\item $\sum_{p\in Z} n_p^2=m^2$
\item For each $L \in \A_w,$ $\sum_{p \in Z \cap L} n_p=m.$
\end{enumerate}
\end{lem}

\subsection{Determinantal syzygies and factoring divisors}
One simple way in which linear syzygies can arise comes from a factorization
of divisors. First, a definition
\begin{defn}
A matrix of linear forms is $1-generic$ if it has no zero entry, and
cannot be transformed by row and column operations to have a zero entry.
\end{defn}
For $Y \subseteq \PP^n$ irreducible and linearly normal, if 
there exist line bundles ${\mathcal L}_1$ and ${\mathcal L}_2$ such that 
$\mathcal{O}_Y(1) = {\mathcal L}_1 \otimes {\mathcal L}_2$ with $h^0({\mathcal L}_i)= a_i$, then the $a_1 \times a_2$ matrix $\gamma$
representing the multiplication table 
\[
H^0({\mathcal L}_1)\otimes H^0({\mathcal L}_2) \longrightarrow H^0(\mathcal{O}_Y (1)) 
\] 
is $1$--generic. More explicitly (see \cite{eis}), if 
\[
H^0({\mathcal L}_1) = Span_{\C}\{e_1,\ldots,e_{a_1}\} \mbox{ and } H^0({\mathcal L}_2) = Span_{\C}\{f_1,\ldots,f_{a_2}\},
\]
then $\gamma$ has $(i,j)$ entry $e_i \otimes f_j$, corresponding 
to a linear form on $\PP^n$, and elements of the ideal $I_2(\gamma)$ of 
$2 \times 2$ minors of $\gamma$ vanish on $Y$. 
The most familiar example occurs when $a_1=2$ and $a_2=k$. In this case,
the minimal free resolution of $I_2(\gamma)$ is an Eagon-Northcott complex. 
This relates to geometry via scrolls: let $\Psi$ be the locus of points 
where a $1$--generic matrix
\[
\gamma =\left[ \!
\begin{array}{cccc}
l_1 & \cdots & l_k\\
m_1 & \cdots & m_k\
\end{array}\! \right]
\]
has rank one. If 
\[
L_{[\lambda:\nu]} = \{p \in \PP^n \mid \lambda l_1(p)+\mu m_1(p) = \cdots = 
\lambda l_k(p)+\mu m_k(p) = 0 \},
\]
then (see 9.10 of \cite{harris})
\[
\Psi = \bigcup\limits_{[\lambda:\nu]\in \PP^1} L_{[\lambda:\nu]},
\] 
where $ L_{[\lambda:\nu]}\simeq \PP^{n-k}$, so $\Psi$ is a union of linear spaces.
Geometrically, the zero locus of the $2 \times 2$ minors of $\gamma$ is 
a scroll which contains $V(I_Y)$.
\subsection{Connecting nets and determinantal syzygies}
The computation in the proof
of Theorem~\ref{sectionsD} and the fact that $h^1(D) \ge 0$ shows that
if $D_{\A} = A+B$ with $A = mE_0 - \sum a_iE_i$, then
\[
h^0(A) \ge {m+2 \choose 2} -\!\!\sum\limits_{p \in L_2(\A)} \!\!\!{a_i+1 \choose 2},\mbox{ } 
h^0(B) \ge {d+1-m \choose 2} -\!\!\sum\limits_{p \in L_2(\A)} \!\!\!{\mu(p) - a_i +1 \choose 2}.
\]
For an arrangement $\A$, if there exists a choice of parameters $m$ and $a_i$ such that
$h^0(A) = a \ge 2$ and $h^0(B)=b \ge 3$, then the results of the previous section show 
that there will exist linear first syzygies on $I$.
\begin{exm}\label{NetSections}
We revisit Example~\ref{exm:p1}. Let $A = 3E_0 - \sum_{\{p | \mu(p) = 2\}} E_p$. 
Clearly $A^2 = AK = 0$, so we can only guarantee that $h^0(A) \ge 1$. In
fact, $h^0(A)=2$, hence $h^1(A) = 1$. To see this, note that a direct 
computation shows that the space of cubics passing through the
nine multiple points of $\A$ is two dimensional. Since
\[
Span_{\C} \Big( L_1L_6L_9, L_3L_4L_7, L_2L_5L_8 \Big) \subseteq H^0(A),
\]
and any two of these are independent, we see that the sections are given
by the net. Next, consider the residual divisor $B = D_{\A}-A$. Since 
$B^2 = 16-18 = -2$ and $-BK = 15 - 9 = 6$, we have that $h^0(B) \ge 3$. 
In fact, equality holds, so $I$ contains the $2 \times 2$ 
minors of a $2 \times 3$ matrix of linear forms, explaining the linear syzygies. 
\end{exm}
\begin{lem}\label{weakMN}
A $(k,m)$ multinet gives a divisor $A$ on $X$ such that $h^0(A)=2$.
\end{lem}
\begin{proof}
Let 
\[
A = mE_0 - \sum\limits_{p \in Z}n_pE_p.
\]
Condition (1) of Definition~\ref{MNet} implies that for each block $A_i$ of the multinet, 
$\prod_{L \in A_i} L^{w(L)}$ is homogeneous of degree $m$, 
and Condition (3) shows that it vanishes to order exactly $n_p$ on $E_p$. 
In particular, this shows that $\prod_{L \in A_i} L^{w(L)} \in H^0(A)$. If 
all the blocks of the partition were independent, then this would imply
that $h^0(A) \ge k$, but it turns out that the sections are all fibers of
a pencil of plane curves, which follows from Theorem 3.11 of \cite{FY}.
\end{proof}

In \cite{FY}, Falk and Yuzvinsky show that the following are equivalent:
\begin{enumerate}
\item $R^1(A)$ contains a nonlocal component $\simeq \PP^{k-2}$.
\item $\A$ supports a $(k,m)$ multinet.
\item $\exists$ a pencil of plane curves with connected fibers, with at least three
fibers (loci of) products of linear forms, and $\A$ is the union of all such fibers.
\end{enumerate}
In general, determining the dimension of $h^1(D)$ for $D \in Pic(X)$ is not easy. However,
in the special case of a net, there is enough information to give a lower bound for the
dimension of the sections of the residual divisor which is often exact. 
\begin{thm}\label{netSyz}
If $\A$ is a $(k,m)$ net, then $D_{\A} = A+B$, with 
\[
h^0(A) = 2 \mbox{ and }h^0(B) \ge km - {m+1 \choose 2}.
\]
\end{thm}
\begin{proof}
For a $(k,m)$ net, all lines occur with multiplicity one. Let 
\[
A= mE_0 - \sum\limits_{p \in Z}E_p.
\]
By Lemma~\ref{weakMN}, $h^0(A) =2$. Since
$B = D_{\A} - mE_0 + \sum\limits_{p \in Z}E_p$, $\frac{B^2-BK}{2}+1 $
\[
\begin{array}{cc}
=& \frac{(d-m-1)(d-m-1+3) +2}{2} + (\!\!\sum\limits_{p \in L_2(\A)} \!\!\!\mu(p)E_p 
+ \sum\limits_{p \in Z}E_p)(\!\!\sum\limits_{p \in L_2(\A)} \!\!\!(\mu(p)-1)E_p 
+ \sum\limits_{p \in Z}E_p)\\ 

  = & {d+1-m \choose 2} -{d \choose 2} + \sum\limits_{p \in Z}\mu(p)\\
 = & {m \choose 2} -d(m-1) + \sum\limits_{p \in Z}\mu(p).
\end{array}
\]
We now compute that 
\[
\begin{array}{ccc}
\sum\limits_{p \in Z}\mu(p) + |Z| &=& \sum\limits_{p \in Z}(\mu(p)+1)\\
    & = & k \sum\limits_{p \in Z}n_p \\
    & = &  km^2.
\end{array}
\]
The second line follows since $n_p$ lines from each block $A_i$ pass through
$p$, and there are $k$ blocks. The third line follows from 
Lemma~\ref{numerology} and the fact that $n_p =1$ for a net, 
hence $n_p^2 = n_p$. Since for a $(k,m)$--net $|Z|=m^2$, 
\[
\sum\limits_{p \in Z}\mu(p) = (k-1)m^2 = dm - m^2.
\]
Combining this with the previous calculation shows that for a $(k,m)$--net
\[
h^0(B) = h^1(B) +  {m \choose 2} -d(m-1) +  dm - m^2 \ge d - {m+1 \choose 2}.
\]
Recalling that $km=d$ concludes the proof.
\end{proof}
\begin{cor}
If $\A$ is a $(k,m)$ net with $k \ge m$, then $I$ contains the 
$2 \times 2$ minors of a $1$-generic $2 \times \Big( km - {m+1 \choose 2} \Big)$ matrix. 
Thus the resolution of $I$ contains an Eagon-Northcott complex as a subcomplex.
\end{cor}
\begin{proof}
Since $k \in \{3,4\}$, if $(k,m) = (3,2)$ or $(3,3)$ then by Theorem~\ref{netSyz}
$h^0(B) \ge 3$, and if $(k,m) = (4,3)$ or $(4,4)$ then $h^0(B) \ge 6$. Note that 
the only known example of a $4$-net is the $(4,3)$ net corresponding to
the Hessian configuration.
\end{proof}
\begin{exm}\label{ShowMeThm}
For the arrangement $A_3$ appearing in Example~\ref{exm:braidex}, $Z$ is the 
collection of multiple points, and  
\[
A = 2E_0-\!\!\sum\limits_{\{p | \mu(p) = 2\}} \!\!E_p
\]
and
\[
B = 3E_0-\!\!\sum\limits_{p \in L_2(\A)} \!\!E_p.
\]
So $d - {m+1 \choose 2} = 6-3 =3$ and $I$ 
contains the $2\times 2$ minors of a $2 \times 3$ matrix. 
\end{exm}
\section{Connection to Derivations}\label{sec:five}
In this section, we show that the generators of the Jacobian ideal
of $\A \subseteq \PP^2$ are contained in $H^0(D_{\A})$, and that the
associated projection map $X\rightarrow \PP^2$ has degree 
$\sum_{p \in L_2(\A)} \mu(p) -|\A| +1$. This relates $X$ to
one of the fundamental objects in arrangement theory:
the module $D(\A)$ of derivations tangent to $\A$. 
\begin{defn}
 $D({\mathcal A}) = \{ \theta \mid \theta(\alpha_i) \in \langle \alpha_i \rangle \mbox{ for all }\alpha_i \mbox{  such that }V(\alpha_i) \in \A \}.$
\end{defn}
The module $D(\A)$ is a graded $S = \C[x_0,\ldots,x_n]$ module, 
and over a field of characteristic zero,
$D({\mathcal A})\simeq E \oplus D_0({\mathcal A})$, where
$E$ is the Euler derivation and $D_0({\mathcal A})$
corresponds to the module of syzygies on the Jacobian ideal
$J_{\alpha}$ of the defining polynomial $\alpha = \prod_{i=1}^d \alpha_i$ 
of ${\mathcal A}$. 
An arrangement $\mathcal{A}$ is {\em free} if $D({\mathcal A}) \simeq
\oplus S(-a_i)$; the $a_i$ are called the {\em exponents} of ${\mathcal A}$.
Terao's theorem \cite{t} is that if $D({\mathcal A}) \simeq \oplus S(-a_i)$, 
then $P(M,t)= \prod(1+a_it)$. Supersolvable arrangements are free, so 
Example~\ref{exm:braidex} is of this type, and 
\[
P(A_3,t) = (1+t)(1+2t)(1+3t)
\]
On the other hand, the arrangement $\A$ of 
Example~\ref{exm:secondex} is not free, and $P(\A,t)$
does not factor. However, it is shown in \cite{sAdv} that
for $\A \subseteq \PP^2$ the Poincare polynomial is $(1+t)\cdot c_t(D_0^\vee)$,
where $c_t$ is the Chern polynomial and $D_0^\vee$ is the dual of the
rank two vector bundle associated to $D_0$. An easy localization 
argument \cite{sCMH} shows that in this case, the Jacobian ideal 
is a local complete intersection with
\[
deg(J_{\alpha}) = \!\!\sum\limits_{p \in L_2(\A)} \mu(p)^2.
\]
For $\A \subseteq \PP^n$ with $n\ge 3$, some generalizations are possible, see \cite{ms}. 
\begin{prop}\label{Jac}
For an arrangement $\A \subseteq \PP^2$, 
\[
J_{\alpha} \subseteq L = \langle l_1, \ldots, l_d \rangle =\langle \frac{\alpha}{\alpha_1}, \ldots, \frac{\alpha}{\alpha_d} \rangle.
\]
\end{prop}
\begin{proof}
By Lemma~\ref{idealT}, 
\[
L = \bigcap\limits_{p \in L_2(\A)} I_{p}^{\mu(p)} \subseteq \C[x,y,z].
\]
The ideal $J_{\alpha}$ is generated in degree $d-1$. The result of \cite{sCMH} mentioned
above implies that at any point $p \in L_2(\A)$, the localization 
$(J_{\alpha})_p$ is a local complete intersection: changing coordinates so
that $p = (0:0:1)$, and writing $\A = V(L_0L_1)$ with $L_0$ the product of the defining 
linear forms which vanish at $p$ and $L_1$ the product of the remaining forms, we have 
\[
(J_{\alpha})_p = \langle \partial(L_0)/\partial_x,\partial(L_0)/\partial_y \rangle. 
\]
In particular, both generators are of degree $\mu(p)$, so in the primary
decomposition of $J_{\alpha}$, the primary component associated to $I(p)$ is 
contained in $I_p^{\mu(p)}$. Now,
\[
Sing(\A) = \bigcup\limits_{I(p) \in Ass(\sqrt{J_{\alpha}})} V(I(p)).
\]
If $Q_p$ is the $I(p)$-primary component of $J_{\alpha}$, then 
\[
\bigcap\limits_{I(p) \in S} Q_p \subseteq L.
\]
Since the saturation of $J_{\alpha}$ with respect to $\langle x,y,z \rangle$ is the
left hand side, and $J_{\alpha}$ is generated in degree $d-1$, the result follows.
\end{proof}
\noindent The inclusion $W = J_{\alpha} \subseteq H^0(D_{\A})$ corresponds to
an induced map
\[
\xymatrix{
\PP^2 \setminus V(\alpha) \ar[r]^{\phi_{\A}} \ar[dr]^{\psi}  & \PP(H^0(D_{\A}))\ar[d]^{\pi}\\
                 & \PP(W) }
\]
\begin{prop}\label{projection}
The degree of $\pi$ is $\!\!\sum\limits_{p \in L_2(\A)} \mu(p) -|\A| +1$.
\end{prop}
\begin{proof}
In \cite{DimP}, Dimca and Papadima show that on a projective hyperplane complement 
$\PP^n \setminus V(\alpha)$, the degree of the gradient map $\psi$
\[
\PP^n \setminus V(\alpha) \xrightarrow{\left[ \!\begin{array}{ccc}
\partial(\alpha)/\partial x_0 : & \cdots& :\partial(\alpha)/\partial x_n
\end{array}\! \right]}
\PP^n
\]
is equal to $b_n(\PP^n \setminus V(\alpha))$; for a configuration of 
$\A \subseteq \PP^2$ this means the degree of the gradient map is 
\[
\!\!\sum\limits_{p \in L_2(\A)} \mu(p) -|\A| +1.
\]
By Theorem~\ref{geometrySurf}, $\phi_{A}$ is an isomorphism on $\PP^2 \setminus V(\alpha)$, and the
result follows.
\end{proof}
\vskip .1in
\noindent{\bf Concluding Remarks and Questions} 
\begin{enumerate}
\item To study $Tor^R_i(C(A),\C)_{i+1}$, it suffices to 
restrict to the case of line arrangements. This follows since the quadratic 
generators of $C(A)$ depend only on $L_2(\A)$, hence taking the intersection
of $\A$ with a generic $\PP^2$ leaves these generators (and relations among them) unchanged. 
The graded betti numbers for $I_2$ are not combinatorial invariants; it would
be interesting to understand how the geometry of $\A$ governs $b_{ij}$, even
for $j = i+1$. 

\item Can freeness of $D(\A)$ be related to the surface $X$ and divisor $D_{\A}$? It seems possible that there is a connection between $D_{\A}$ and
multiarrangements, studied recently in \cite{atw},
\cite{yo1}, \cite{yo2}. 

\item For $n \ge 3$, is $C(\A)$ the homogeneous coordinate ring
of a blowup of $\PP^n$ along a locus related to the arrangement $\A$?
Since $C(\A)$ is Cohen-Macaulay, if this is true, then 
combining Riemann-Roch with Terao's result would yield a formula
for the global sections of $D_{\A}$. 
\end{enumerate} 
\vskip .1in
\noindent{\bf Acknowledgements} Computations were performed using Macaulay2,
by Grayson and Stillman, and available at: {\tt http://www.math.uiuc.edu/Macaulay2/}.
Many of the examples in this paper appear in Suciu's 
survey \cite{Su}; that paper is also an excellent reference for the many 
questions on resonance varieties not treated here. This work began during
a visit to the Fields institute for a conference in honour of Peter Orlik.
\bibliographystyle{amsalpha}

\end{document}